\documentclass[12pt]{amsart}

\usepackage [T1] {fontenc}
\usepackage {amssymb}
\usepackage {amsmath}
\usepackage {amsthm}
\usepackage{tikz}
\usepackage{hyperref}
\usepackage{url}

 \usepackage[titletoc,toc,title]{appendix}

\usetikzlibrary{arrows,decorations.pathmorphing,backgrounds,positioning,fit,petri}

\usetikzlibrary{matrix}

\usepackage{enumitem}

\DeclareMathOperator {\pr} {pr}

\DeclareMathOperator {\SL} {SL}

\DeclareMathOperator {\GL} {GL}

\newcommand {\gm} {\ensuremath{\mathbb{G}_{\mathrm{m}}}}
\DeclareMathOperator {\ws} {ws}
\DeclareMathOperator {\zcl} {Zcl}

\newcommand {\C} {\mathbb{C}}

\newcommand {\s} {\mathfrak{S}}

\newcommand {\rar} {\rightarrow}
\newcommand {\seq} {\subseteq}

\theoremstyle {definition}
\newtheorem {definition}{Definition} [section]

\newtheorem* {notation} {Notation}

\theoremstyle {plain}

\newtheorem {lemma} [definition] {Lemma}
\newtheorem {theorem} [definition] {Theorem}
\newtheorem {proposition} [definition] {Proposition}

\newtheorem {conjecture} [definition] {Conjecture}

\theoremstyle {remark}
\newtheorem {remark} [definition] {Remark}

\makeatletter
\newcommand {\forksym} {\raise0.2ex\hbox{\ooalign{\hidewidth$\vert$\hidewidth\cr\raise-0.9ex\hbox{$\smile$}}}}
\def\@forksym@#1#2{\mathrel{\mathop{\forksym}\displaylimits_{#2}}}
\def\forkind{\@ifnextchar_{\@forksym@}{\forksym}}

\makeatother

\calclayout

\begin {document}

\title{Some remarks on atypical intersections}

\author{Vahagn Aslanyan}
\address{Department of Mathematical Sciences, Carnegie Mellon University, 5000 Forbes Ave., Pittsburgh, PA 15213, USA}

\address{\texttt{Current address:} School of Mathematics, University of East Anglia, Norwich, NR4 7TJ, UK}
\email{V.Aslanyan@uea.ac.uk}


\thanks{This work was done while I was a postdoctoral associate at Carnegie Mellon University. Some revisions were made at the University of East Anglia. Partially supported by EPSRC grant EP/S017313/1.}

\date{\today}

\keywords{Unlikely intersection, Zilber-Pink, semiabelian variety, $j$-function, special variety, Ax-Schanuel}

\subjclass[2010]{11G10, 11G18}


\maketitle

\begin{abstract}
In this paper we show how some known weak forms of the Zilber--Pink conjecture can be strengthened by combining them with the Mordell--Lang conjecture or its variants.
We illustrate this idea by proving some theorems on atypical intersections in the semiabelian and modular settings. Given a ``finitely generated'' set $\Gamma$ with a certain structure, we consider $\Gamma$-special subvarieties---weakly special subvarieties containing a point of $\Gamma$---and show that every variety $V$ contains only finitely many maximal $\Gamma$-atypical subvarieties, i.e. atypical intersections of $V$ with $\Gamma$-special varieties the weakly special closures of which are $\Gamma$-special.
\end{abstract}



\section{Introduction}

\subsection{The Zilber--Pink conjecture} The Zilber--Pink conjecture is a statement about \emph{atypical} intersections of an algebraic variety with some (countable) collection of \emph{special} varieties. An intersection is \emph{atypical} or \emph{unlikely} if its dimension is larger than expected. 
The Zilber--Pink conjecture states, roughly, that atypical intersections of a variety with special varieties are governed by finitely many special varieties (precise definitions and statements will be given shortly).

The conjecture for algebraic tori and, more generally, for semiabelian varieties was first posed by Zilber in his work on Schanuel's conjecture and the model theory of complex exponentiation \cite{Zilb-exp-sum-published}. He showed, in particular, that it implies the Mordell--Lang conjecture. Bombieri, Masser and Zannier \cite{Bom-Mas-Zan} gave an equivalent formulation independently. Pink \cite{Pink,Pink-2} proposed (again independently) a more general conjecture for mixed Shimura varieties which also implies the Andr\'e--Oort conjecture. 


Let us start with a rigorous definition of atypical intersections. 
Let $V$ and $W$ be subvarieties of some variety $S$. A non-empty component $X$ of the intersection $V \cap W$ is \emph{atypical} in $S$ if $\dim X > \dim V + \dim W - \dim S$, and \emph{typical} if $\dim X = \dim V + \dim W - \dim S$. Note that if $S$ is smooth then a non-strict inequality always holds.


Now let us describe \emph{special} varieties. For a semiabelian variety $\mathfrak{S}$ (defined over $\mathbb{C}$) its \emph{special} subvarieties are torsion cosets of semiabelian subvarieties of $\mathfrak{S}$, and arbitrary cosets are called \emph{weakly special} subvarieties. Note that special subvarieties are precisely the irreducible components of algebraic subgroups of $\mathfrak{S}$. In the modular setting, the \emph{special} subvarieties of $Y(1)^n$ (where the modular curve $Y(1)$ is identified with the affine line $\mathbb{C}$) are irreducible components of algebraic varieties defined by modular equations, that is, equations of the form $\Phi_N(x_i,x_k)=0$ for some $1\leq i \leq k \leq n$ where $\Phi_N(X,Y)$ is a modular polynomial (see \cite{Lang-elliptic}). If we also allow equations of the form $x_i = c_i$ for constants $c_i \in \mathbb{C}$ then we get \emph{weakly special} subvarieties.


Now let $\mathfrak{S}$ be a semiabelian variety or $Y(1)^n$, and let $S$ be a special subvariety of $\mathfrak{S}$. For a subvariety $V \subseteq S$ an \emph{atypical subvariety} of $V$ in $S$ is an atypical (in $S$) component $X$ of an intersection $V\cap T$ where $T \subseteq S$ is special. When we do not specify $S$ then we mean $S = \mathfrak{S}$, i.e. an atypical subvariety of $V$ is an atypical subvariety of $V$ in $\mathfrak{S}$.

Now we are ready to formulate the Zilber--Pink conjecture for $\mathfrak{S}$. There are many equivalent forms of the conjecture; we consider two of them (see \cite{Zilb-exp-sum-published,Bom-Mas-Zan,Pink,Habegger-Pila-o-min-certain}).

\begin{conjecture}[Zilber--Pink for $\mathfrak{S}$: Formulation 1]\label{intro-ZP-1}
Let $\mathfrak{S}$ be a semiabelian variety or $Y(1)^n$ and $V \subseteq \mathfrak{S}$ be an algebraic subvariety. Then $V$ contains only finitely many maximal atypical subvarieties.
\end{conjecture}

\begin{conjecture}[Zilber--Pink for $\mathfrak{S}$: Formulation 2]\label{intro-ZP-2}
Let $\mathfrak{S}$ be a semiabelian variety or $Y(1)^n$ and $V \subseteq \mathfrak{S}$ be an algebraic subvariety. Then there is a finite collection $\Sigma$ of proper special subvarieties of $\mathfrak{S}$ such that every atypical subvariety $X$ of $V$ is contained in some $T \in \Sigma$.
\end{conjecture}



Although the Zilber--Pink conjecture is wide open, many special cases and weak versions have been proven in the past two decades. The reader is referred to \cite{Zannier-book-unlikely,pila-ellipt-mod-surf,Habegger-Pila-o-min-certain,Tsimerman-AO-A_g,Daw-Ren,Aslanyan-weakMZPD} for various results and recent developments around this conjecture. We now formulate two well known weak Zilber--Pink theorems which play a key role in this paper.

\begin{theorem}[Weak Zilber--Pink for semiabelian varieties, \cite{Zilb-exp-sum-published,Kirby-semiab,Bom-Mas-Zan}]\label{Weak-CIT-intro}
Let $\mathfrak{S}$ be a semiabelian variety and $V$ be an algebraic subvariety of $\mathfrak{S}$. Then atypical components of intersections of $V$ with cosets of algebraic subgroups of $\mathfrak{S}$ are contained in cosets of finitely many algebraic subgroups.
\end{theorem}


\begin{theorem}[Weak Modular Zilber--Pink, \cite{Pila-Tsim-Ax-j,Aslanyan-weakMZPD}]\label{weak-modularZP-intro}
Every algebraic subvariety $V \subseteq Y(1)^n$ contains only finitely many maximal strongly atypical subvarieties, that is, atypical subvarieties with no constant coordinate.
\end{theorem}

\subsection{$\Gamma$-special varieties and the Mordell--Lang conjecture}





\begin{definition}
\begin{itemize}[leftmargin=0.5cm]
    \item[] 
    \item A subset of a semiabelian variety is said to be \emph{a structure of finite rank} if it is a subgroup of finite rank, that is, $\dim_{\mathbb{Q}} (\Gamma \otimes \mathbb{Q})$ is finite.
    \item For a set $\Xi\seq \C$ we denote by $\overline{\Xi}$ the union of all Hecke orbits of points of $\Xi$, that is, $\overline{\Xi} = \{ \eta\in \C: \Phi_N(\xi,\eta)=0 \mbox{ for some } N\in \mathbb{N},~ \xi\in \Xi \}$. A subset $\Gamma \seq Y(1)^n(\C)$ is called \emph{a structure of finite rank} if there is a set $\Xi\seq \C$ containing only finitely many non-special points such that $\Gamma = \left(\overline{\Xi}\right)^n.$

\end{itemize}
\end{definition}

\begin{definition}
Let $\mathfrak{S}$ be a semiabelian variety or $Y(1)^n$ and let $\Gamma \subseteq \mathfrak{S}$ be a structure of finite rank. 
\begin{itemize}[leftmargin=0.5cm]

     \item For an irreducible subvariety $X \subseteq \mathfrak{S}$, \emph{the weakly special closure} of $X$, denoted $\langle X \rangle_{\ws}$, is the smallest weakly special subvariety containing $X$. Similarly, $\langle X \rangle$ denotes \emph{the special closure} of $X$, i.e. the smallest special subvariety containing $X$.
     
    \item A weakly special subvariety of $\s$ is called $\Gamma$-\emph{special} if it contains a point of $\Gamma$.
    
     \item Given varieties $V \seq S \seq \s$, with $S$ special, a \emph{(weakly) atypical} subvariety of $V$ in $S$ is an atypical component (in $S$) of an intersection of $V$ with a (weakly) special subvariety of $S$.
    
    \item A weakly atypical subvariety $X \seq V$ is $\Gamma$-atypical if $\langle X \rangle_{\ws}$ is $\Gamma$-special.
\end{itemize}
\end{definition}

In terms of $\Gamma$-special varieties the Mordell--Lang conjecture can be stated as follows.

\begin{theorem}[Mordell--Lang for $\s$]\label{mordell-lang-combined}
Let $\mathfrak{S}$ be a semiabelian variety or $Y(1)^n$ and let $\Gamma \subseteq S$ be a structure of finite rank. Then every algebraic variety $V \subseteq \mathfrak{S}$ contains only finitely many maximal $\Gamma$-special subvarieties.
\end{theorem}

For semiabelian varieties this was proven in a series of papers by Faltings, Vojta, McQuillan and many others (see \cite{mcquillan}). Its modular analogue was established by Pila (see \cite[Theorem 6.6]{pila-ellipt-mod-surf}) generalising an earlier result of Habegger and Pila from \cite{Habegger-Pila-beyond}.

\subsection{Main results and key ideas of the proofs}

In this paper we demonstrate how the aforementioned weak Zilber--Pink type theorems can be combined with Mordell--Lang to generalise the former and thus establish new variants of the Zilber--Pink conjecture. We prove some precise results (stated below) in the semiabelian and modular settings to illustrate this idea. Moreover, we believe our methods can be extended to work in the more general context of Shimura varieties. Nevertheless, we choose to work with semiabelian varieties and products of modular curves to keep the paper short and simple. 

Our main results can be combined into the following theorem.

\begin{theorem}\label{gamma-cit-intro}
Let $\mathfrak{S}$ be a semiabelian variety or $Y(1)^n$, let $\Gamma\subseteq \mathfrak{S}$ be a structure of finite rank and let $S\subseteq \mathfrak{S}$ be a $\Gamma$-special subvariety. Then every subvariety $V \subseteq S$ contains only finitely many maximal $\Gamma$-atypical subvarieties in $S$.
\end{theorem}

This theorem generalises Theorems \ref{Weak-CIT-intro} and \ref{weak-modularZP-intro}. Observe that the latter states that $V$ contains finitely many maximal atypical subvarieties with no constant coordinates, and Theorem \ref{gamma-cit-intro} for $Y(1)^n$ shows that we can also deal with atypical subvarieties with constant coordinates provided that we limit those constants to a small set. In particular, $V$ contains only finitely many maximal atypical subvarieties all constant coordinates of which are special. In terms of optimal varieties (see Section \ref{section-optimal}) this is equivalent to the statement that $V$ contains only finitely many optimal subvarieties whose weakly special closures are special. This statement generalises \cite[Corollary 9.11]{Habegger-Pila-o-min-certain}. 

Note that Pila and Scanlon have recently proven some differential algebraic Zilber--Pink theorems where they work over a differential field $(K; D)$ and consider atypical intersections possibly with constant coordinates which are not constant in the differential algebraic sense, i.e. they allow equations $x_i = c_i$ where $c_i \in K$ with $Dc_i \neq 0$. In particular, $c_i$ cannot be algebraic (over $\mathbb{Q}$) since algebraic numbers are constant in any differential field. See Scanlon's slides \cite{Pila-Scanlon-dif-ZP} for details.

Let us outline the strategy of the proof of Theorem \ref{gamma-cit-intro} assuming for simplicity that $S=\mathfrak{S}$ is a semiabelian variety. Given a subvariety $V$ and an algebraic subgroup $T$ of $\mathfrak{S}$, we show that a generic coset of $T$ intersects $V$ typically (or does not intersect it at all). This is consistent with the intuitive idea that ``generic'' varieties intersect typically. Thus, the set of all cosets $U$ for which $V \cap U$ is atypical in $\mathfrak{S}$ is a constructible subset $C_T$ of the quotient $\mathfrak{S}/T$ of lower dimension. If we restrict to $\Gamma$-atypical subvarieties then we can use the Mordell--Lang conjecture to deduce that $C_T \cap \Gamma'$ is contained in the union of finitely many $\Gamma'$-special subvarieties of $C_T$ where $\Gamma'$ is the image of $\Gamma$ in $\s/T$ under the natural projection. On the other hand, by Theorem \ref{Weak-CIT-intro} we need to consider only finitely many subgroups $T$ which yields the desired result. 
We also show how the uniform version of Theorem \ref{gamma-cit-intro} can be deduced from the uniform versions of weak Zilber--Pink and Mordell--Lang.

Note that our arguments are quite general and should go through in other settings too provided there is an Ax--Schanuel theorem (which is the key ingredient in the proofs of Theorems \ref{Weak-CIT-intro} and \ref{weak-modularZP-intro}) and some analogue of the Mordell--Lang or Andr\'e--Oort conjectures. Furthermore, Daw and Ren showed in \cite{Daw-Ren} that the Zilber--Pink conjecture for Shimura varieties can be reduced to a conjecture on finiteness of optimal points. It seems their methods can be adapted to reduce Theorem \ref{gamma-cit-intro} (at least for $Y(1)^n$) to a similar point counting problem which would follow from Mordell--Lang, and that will then give another proof of that theorem. Daw has shown in a private communication to me that this can indeed be done when $\Gamma$ is the set of special points. See Section \ref{section-optimal} for more details.

\section{$\Gamma$-atypical subvarieties in semiabelian varieties}\label{section-semiab}

In this section $\s$ denotes a semiabelian variety, written additively. However, algebraic tori are written multiplicatively since they are subgroups of a multiplicative group $\gm^n$.

The following simple fact (and its obvious analogue in the modular setting) will be used repeatedly in the paper.

\begin{lemma}\label{lem-ws-atyp}
Let $V\subseteq \mathfrak{S}$ be a subvariety. If $X$ is a weakly atypical subvariety of $V$ in $\mathfrak{S}$ then $X$ is an atypical component of the intersection $V \cap \langle X \rangle_{\ws}$ in $\mathfrak{S}$.
\end{lemma}
\begin{proof}
Assume $T \subseteq \mathfrak{S}$ is weakly special such that $X$ is an atypical component of $V \cap T$ in $\mathfrak{S}$. Then $\langle X \rangle_{\ws}\subseteq T$ and so
$$\dim X > \dim V + \dim T - \dim \mathfrak{S} \geq \dim V + \dim \langle X \rangle_{\ws} - \dim \mathfrak{S}.$$
Now if $Y \subseteq V \cap \langle X \rangle_{\ws}$ is a component containing $X$ then $Y \subseteq V \cap T$. Since $X$ is an irreducible component of $V\cap T$, so is $Y$ and in fact $X=Y$.
\end{proof}

The analogous statement for atypical subvarieties and special closures holds too.

Let $T \subseteq \mathfrak{S}$ be an algebraic subgroup and $V \subseteq \mathfrak{S}$ be an irreducible algebraic subvariety. We show that generic cosets of $T$ intersect $V$ typically. First, note that the quotient $\mathfrak{S}/T$ is (definably isomorphic to) an algebraic group and the natural projection $\pi : \mathfrak{S} \rightarrow \mathfrak{S}/T$ is a morphism of algebraic groups.\footnote{Note that this follows from elimination of imaginaries in algebraically closed fields and the fact that constructible groups are definably isomorphic to algebraic groups. See \cite[Chapter 7]{Mar}.} Moreover, $\mathfrak{S}/T$ is connected and hence irreducible. 

\begin{lemma}\label{cit-generic-typical}
Let $T$ and $V$ be as above. The set \[ C := C_T := C_{T,V} := \{ u \in \mathfrak{S}/T: V \cap \pi^{-1}(u) \mbox{ is atypical in } \mathfrak{S} \} \] is constructible and not Zariski dense in $\mathfrak{S}/T$.
\end{lemma}
Note that by definition, atypicality of an intersection implies that it is non-empty, hence if $V \cap \pi^{-1}(u) = \emptyset$ then $u \notin C$.
\begin{proof}
Let $\theta:=\pi|_V:V \rar \s/T$ be the restriction of of $\pi$ to $V$. Observe that for every $u \in \mathfrak{S}/T$ we have $\dim \pi^{-1}(u) = \dim T$, for $\pi^{-1}(u)$ is a coset of $T$. Hence \[ C = \{ u \in \mathfrak{S}/T : \dim \theta^{-1}(u) > \dim V + \dim T - \dim \mathfrak{S}  \}\] which is constructible since $\theta$ is a morphism of varieties and $\dim V + \dim T - \dim \mathfrak{S}$ is a fixed number independent of $u$. 

Now assume $C$ is Zariski dense in $\s/T$. Pick a generic point $w \in \mathfrak{S}/T$. Then $w \in C$, and so $V\cap \pi^{-1}(w) \neq \emptyset$, hence $w\in \pi(V)$. This means $\theta$ is a dominant map as its image contains a generic point of $\s/T$. Therefore, by the fibre dimension theorem (\cite[Theorem 1.25]{Shafarevich}),
$ \dim \theta^{-1}(w) = \dim V - \dim \s/T = \dim V + \dim T - \dim \s. $
Hence $w \notin C$, which is a contradiction.
\end{proof}

The following is the central theorem of this section which implies some related results.

\begin{theorem}\label{gamma-cit}
Let $\mathfrak{S}$ be a semiabelian variety and let $\Gamma\seq \s$ be a subgroup of finite rank. Then for every subvariety $V \subseteq \mathfrak{S}$ there is a finite collection $\Sigma$ of proper $\Gamma$-special subvarieties of $\mathfrak{S}$ such that any $\Gamma$-atypical subvariety of $V$ (in $\mathfrak{S}$) is contained in some $T \in \Sigma$.
\end{theorem}

\begin{proof}
It is easy to see that an atypical subvariety of $V$ in $\mathfrak{S}$ is also an atypical subvariety of an irreducible component of $V$. Hence we may assume $V$ is irreducible.

Let $\Sigma_0$ be the finite collection of algebraic subgroups of $\mathfrak{S}$ given by \cite[Theorem 4.6]{Kirby-semiab} (which is a stronger version of Theorem \ref{Weak-CIT-intro}) for $V$. Let further $X$ be a $\Gamma$-atypical subvariety of $V$. 
Then $\langle X \rangle_{\ws}$ is $\Gamma$-special and $X$ is an atypical component of $V \cap \langle X \rangle_{\ws}$. 

By \cite[Theorem 4.6]{Kirby-semiab}, there is $b \in \mathfrak{S}$ and $T \in \Sigma_0$ such that $X \subseteq b+T$. Hence $\langle X \rangle_{\ws} \subseteq b+T$ and so $b+T = \gamma + T$ for some $\gamma \in \Gamma$. Further, we also have
\begin{gather*}
    \dim V + \dim \langle X \rangle_{\ws} - \dim \mathfrak{S} < \dim X =\\
     \dim (V\cap (\gamma +T)) + \dim (\langle X \rangle_{\ws} \cap (\gamma +T)) - \dim (\gamma +T) =\\ 
     \dim (V\cap (\gamma +T)) + \dim \langle X \rangle_{\ws} - \dim (\gamma+T).
\end{gather*}
Thus, $V$ and $\gamma + T$ intersect atypically in $\mathfrak{S}$. Hence $\pi_T(\gamma)\in C_T$ where $C_T$ is defined as in Lemma \ref{cit-generic-typical} and $\pi_T: \s\rar \s/T$ is the natural projection. 

Now we apply the Mordell--Lang theorem to the Zariski closure $C_T^{\zcl}$ of $C_T$ and the finite rank group $\pi_T(\Gamma)\seq \s/T$. We get a finite collection $\Delta_T$ of maximal $\pi_T(\Gamma)$-special subvarieties of $C_T^{\zcl}$. 
This means that $\pi_T(\gamma)\in A$ for some $A\in \Delta_T$. Hence $X\seq \gamma+T \seq \pi_T^{-1}(A)\subsetneq \s$ and so $X$ is contained in an irreducible component of $\pi_T^{-1}(A)$ which is a proper $\Gamma$-special subvariety of $\s$. Thus, we may choose  $\Sigma$ to be the finite collection of $\Gamma$-special irreducible components of all cosets $\pi_T^{-1}(A)$ for $T \in \Sigma_0$ and $A\in \Delta_T$.
\end{proof}

\begin{theorem}\label{gamma-semiab-special}
Let $\mathfrak{S}$ be a semiabelian variety, let $\Gamma\subseteq \mathfrak{S}$ be a subgroup of finite rank, and let $S \subseteq \mathfrak{S}$ be a $\Gamma$-special subvariety. Then for every subvariety $V \subseteq S$, there is a finite collection $\Sigma$ of proper $\Gamma$-special subvarieties of $S$ such that any $\Gamma$-atypical subvariety of $V$ in $S$ is contained in some $T \in \Sigma$.
\end{theorem}
\begin{proof}(cf. \cite[Theorem 4.6]{Kirby-semiab})
Let $S = \gamma + \mathfrak{S}_0$ where $\mathfrak{S}_0$ is a semiabelian subvariety of $\mathfrak{S}$. If $X$ is an atypical component of $V \cap T$ in $S$, where $T\subseteq S$ is $\Gamma$-special, then $X - \gamma$ is an atypical component of $(V-\gamma) \cap (T - \gamma)$ in $\mathfrak{S}_0$. Set $\Gamma_0 := \Gamma \cap \mathfrak{S}_0$. Then $T-\gamma \subseteq \mathfrak{S}_0$ is $\Gamma_0$-special and $X - \gamma$ is $\Gamma_0$-atypical. Let $\Sigma_0$ be the finite set of $\Gamma_0$-special subvarieties of $\mathfrak{S}_0$ given by Theorem \ref{gamma-cit}. Then we can choose $\Sigma = \{ \gamma + T': T' \in \Sigma_0 \}$.
\end{proof}

\begin{remark}\label{remark-maximal}(\textbf{Proof of Theorem \ref{gamma-cit-intro} for semiabelian varieties})
Let $\mathfrak{S},~ \Gamma,~ S$ and $V$ be as above, and let $\Sigma$ be the finite collection of proper $\Gamma$-special subvarieties of $S$ obtained by Theorem \ref{gamma-semiab-special}. Assume $X\subseteq V$ is a maximal $\Gamma$-atypical subvariety in $S$. Then $X\subseteq T$ for some $T \in \Sigma$, hence there is a component $Y$ of $V \cap T$ with $X \subseteq Y$. If $Y$ is an atypical component of $V \cap T$ in $S$ then $X = Y$. So assume $\dim Y = \dim V + \dim T - \dim S$. On the other hand, $\langle X \rangle_{\ws}\subseteq T$ is $\Gamma$-special and $X$ is an atypical component of $V \cap \langle X \rangle_{\ws}$ in $S$. We claim that $X$ is an atypical component of $Y \cap \langle X \rangle_{\ws}$ in $T$. To this end observe that
$$ \dim Y + \dim \langle X \rangle_{\ws} - \dim T = \dim V + \dim \langle X \rangle_{\ws} - \dim S < \dim X.$$
Since $\dim T < \dim S$, we can proceed by induction on $\dim S$.
\end{remark}

\section{$\Gamma$-atypical subvarieties in $Y(1)^n$}\label{section-modular}

In this section we work in a product of modular curves $\s = Y(1)^n$ and identify it with $\C^n$. 
We introduce a piece of notation before proceeding.

\begin{notation}
Let $n$ be a positive integer.
\begin{itemize}[leftmargin=0.5cm]
    \item We write $(n)$ for $(1,\ldots,n)$. The notation ${i} = (i_1, \ldots, i_m) \subseteq (n)$ means that $1\leq i_1 < \ldots < i_m \leq n$, and ${k} = (k_1, \ldots, k_{n-m}) = (n)\setminus {i}$ is the unique tuple $k \subseteq (n)$ such that $\{ 1, \ldots, n \} = \{ i_1, \ldots, i_m \} \cup \{ k_1, \ldots, k_{n-m} \}$.
    
    \item For ${i} = (i_1, \ldots, i_m)\subseteq (n)$ we define $\pr_{{i}}: \mathbb{C}^n \rightarrow \mathbb{C}^m$ to be the projection map onto the ${i}$-coordinates.

    \item For ${c}\in \mathbb{C}^m$, ${i}=(i_1, \ldots, i_m)\subseteq (n)$ and  $Y\seq \C^n$ set $Y_{{i},{c}} := Y \cap \pr_{{i}}^{-1}({c})\subseteq \mathbb{C}^n$.
\end{itemize}
\end{notation}

\begin{lemma}\label{atypical-closed-modular}
Let $S\subseteq Y(1)^n(\mathbb{C})$ be a weakly special variety and let $V \subseteq S$ be an irreducible algebraic subvariety. Fix ${i} = (i_1, \ldots, i_k) \subseteq (n)$ and set $T := \pr_{{i}} S$. Then
\[ C := C_{{i}} := C_{i,V} := \{ {c}\in T : V \cap S_{{i},{c}} \mbox{ is atypical in } S \}\] is constructible in $T$ and $\dim C < \dim T$.
\end{lemma}

\begin{proof}
Let $\theta:=\pr_i|_V:V \rar T$ be the restriction of $\pr_i$ to $V$. It is easy to see that $\dim S_{{i},{c}} = \dim S - \dim T$ for any $c\in T$. Hence $C = \{ {c}\in T: \dim \theta^{-1}(c) > \dim V  - \dim T \}$ is constructible.

Suppose $C$ is dense in $T$. Then $C$ contains a generic point $b$ of $T$. Clearly, $b\in \theta(V)$ which means $\theta$ is dominant. Since $V$ is irreducible, by the fibre dimension theorem we have
$\dim \theta^{-1}(b) = \dim V - \dim T .$
Therefore ${b}\notin C$. This contradiction shows that $C$ cannot be Zariski dense in $T$.
\end{proof}

\begin{definition}[cf. {\cite[Definition 3.8]{Habegger-Pila-o-min-certain}}]
For a weakly special variety $S$ the largest number $N$ for which $\Phi_N$ occurs in the definition of $S$ is called the \emph{complexity} of $S$ and is denoted by $\Delta(S)$.
\end{definition}

\begin{remark}
For a positive integer $N$ there are only finitely many strongly special varieties of complexity at most $N$.
\end{remark}

\begin{proposition}\label{modular-ZP-weakly-atyp}
Given an algebraic subvariety $V$ of a weakly special variety $S$ in $\mathbb{C}^n$, there is a positive integer $N$ such that for every weakly atypical subvariety $X$ of $V$ there is a proper weakly special subvariety $T$ of $S$ with $\Delta(T) \leq N$ such that $X \subseteq T$ and $V \cap T$ is atypical in $S$.
\end{proposition}
\begin{proof}
If $X$ is strongly atypical then it is contained in one of the finitely many special subvarieties of $S$ given by Theorem \ref{weak-modularZP-intro}. Assume $X$ has some constant coordinates, namely, $x_{i_l} = c_l$ for $l=1, \ldots, m$. Let ${i} := (i_1, \ldots, i_m),~ {c} := (c_1, \ldots, c_m)$. Observe that if a constant coordinate is related by a modular equation to another coordinate on $X$, then the latter must also be constant on $X$ for it is irreducible. Therefore, there is no modular relation between an $i$-coordinate and an $(n)\setminus i$-coordinate on $S$. In particular, $S_{i,c}$ is irreducible and hence weakly special. If $V \cap S_{{i}, {c}}$ is atypical in $S$, and hence $S_{i,c}\subsetneq S$, then we can choose $T = S_{{i}, {c}}$. So assume it is a typical intersection, i.e.
$$\dim (V \cap S_{{i}, {c}}) = \dim V + \dim S_{{i}, {c}} - \dim S.$$
Let ${k} := (n)\setminus {i}$ and define $S' := \pr_{{k}} S$ and $V' := \pr_{{k}} V_{{i}, {c}},~ X' := \pr_{{k}} X_{{i}, {c}}$. Then $S' = \pr_{{k}}S_{{i}, {c}}$. Moreover, $S'$ and $X'$ do not have any constant coordinates and $S'$ is strongly special. If $P:= \langle X \rangle_{\ws}$ is the weakly special closure of $X$ then $X$ is an atypical component of $V \cap P$ in $S$, and $P = P_{{i}, {c}}$. Now if $P' := \pr_{{k}} P$ then we claim that $X'$ is an atypical component of $V' \cap P'$ in $S'$. To this end notice that $\dim X' = \dim X,~ \dim V' = \dim V_{{i}, {c}} = \dim (V \cap S_{{i}, {c}}),~ \dim P' = \dim P,~ \dim S' = \dim S_{{i}, {c}}$. Therefore
\begin{gather*}
    \dim X' = \dim X > \dim V + \dim P - \dim S =\\
    (\dim V' - \dim S' + \dim S) +  \dim P' - \dim S = \dim V' + \dim P' - \dim S'.
\end{gather*}
Since $X'$ does not have constant coordinates, we conclude that it is a strongly atypical subvariety of $V'$ in $S'$. On the other hand, $V'$ is a member of a parametric family of varieties depending only on $V$, hence by \cite[Theorem 5.2]{Aslanyan-weakMZPD} (which is the uniform version of Theorem \ref{weak-modularZP-intro}) there is a natural number $N$, depending only on $V$ and $S$ and independent of $c$, and a special subvariety $T' \subseteq S'$ with $\Delta(T') \leq N$ such that $X' \subseteq T'$ and $V' \cap T'$ is atypical in $S'$. Let $T := \pr_{{k}}^{-1}(T') \cap S_{{i},{c}}$. Then $T \subsetneq S$ is weakly special, $\Delta(T) \leq N$, $X \subseteq T$ and $V \cap T$ is atypical in $S$, for
\begin{gather*}
    \dim (V \cap T) = \dim (V_{{i},{c}} \cap T_{{i},{c}}) = \dim (V' \cap T') > \\
     \dim V' + \dim T' - \dim S' = \dim V + \dim T - \dim S.
\end{gather*}
This finishes the proof.
\end{proof}

Now we can state and prove the main result of this section.

\begin{theorem}\label{gamma-modular-ZP}
Let $\Gamma\seq Y(1)^n(\C)$ be a structure of finite rank and let $S$ be a $\Gamma$-special variety. Then for every subvariety $V \subseteq S$ there is a finite collection $\Sigma$ of proper $\Gamma$-special subvarieties of $S$ such that any $\Gamma$-atypical subvariety of $V$ is contained in some $T \in \Sigma$.
\end{theorem}

\begin{proof}
As in the proof of Theorem \ref{gamma-cit}, we may assume $V$ is irreducible.

Let $X \subseteq V$ be $\Gamma$-atypical. Then its weakly special closure $\langle X \rangle_{\ws}$ is $\Gamma$-special. By Proposition \ref{modular-ZP-weakly-atyp} there is a weakly special $T \subsetneq S$ with $\Delta (T) \leq N$ and $X \subseteq T$ where $N$ depends only on $V$ and $S$. Moreover, $V\cap T$ is atypical in $S$. Since $\langle X \rangle_{\ws} \subseteq T$ and $\langle X \rangle_{\ws}$ contains a $\Gamma$-special point, so does $T$ and hence it is $\Gamma$-special. Assume that $x_{i_l} = \gamma_l,~ l=1, \ldots, m,$ are the constant coordinates of $T$ which are not constant on $S$. Set $i := (i_1,\ldots, i_m)$, ${k} := (n) \setminus {i}$ and  $\tilde{T} := \pr_{{k}}^{-1} (\pr_{{k}} T)$, i.e. $\tilde{T}$ is the special subvariety of $S$ defined by the equations of $T$ apart from the equations $x_{i_l}=\gamma_l$. If $\tilde{T} \subsetneq S$ then $\tilde{T}$ is a proper $\Gamma$-special subvariety of $S$ containing $X$ and $\tilde{T}$ belongs to a finite collection $\Theta$ of $\Gamma$-special subvarieties depending only on $V$ and $S$ since $\Delta(\tilde{T})\leq N$. 

Now assume $\tilde{T} = S$. Then $T = S_{{i},{\gamma}}$ where $\gamma = (\gamma_1,\ldots,\gamma_m)\in \pr_i \Gamma$, and $V \cap S_{{i},{\gamma}}$ is atypical in $S$.  Let $C_{{i}} \subseteq \pr_{{i}}S$ be defined as in Lemma \ref{atypical-closed-modular}. Then $\gamma \in C_{{i}}$. Let $\Xi_{{i}}$ be the finite collection of maximal $\pr_i \Gamma$-special subvarieties of $C_{{i}}^{\zcl}$ given by modular Mordell--Lang. Then there is $Q \in \Xi_{{i}}$ with $\gamma\in Q$. Thus, we can choose $\Sigma$ to be the finite collection of all $\Gamma$-special irreducible components of all varieties from
$ \Theta \cup \{ S\cap \pr_{{i}}^{-1} Q : Q \in \Xi_{{i}},~ {i} \subseteq (n) \}.$
\end{proof}

\begin{remark}
We can deduce Theorem \ref{gamma-cit-intro} for $Y(1)^n$ from Theorem \ref{gamma-modular-ZP} as in Remark \ref{remark-maximal}. 
\end{remark}

\section{Uniform versions}\label{section-uniform}

In this section we establish uniform versions of Theorems \ref{gamma-semiab-special} and \ref{gamma-modular-ZP} using uniform versions of Theorems \ref{Weak-CIT-intro}, \ref{weak-modularZP-intro}, and \ref{mordell-lang-combined}. 
In order to combine the statements in the semiabelian and modular settings into one theorem, we introduce a piece of notation.

\begin{definition}
Given two weakly special varieties $S, T\seq Y(1)^n(\C)$ and a set $\Gamma\seq Y(1)^n(\C)$, we say that $T$ is a $\Gamma$-\emph{translate} of $S$ if for some $i\seq (n)$ and for some $\gamma\in \pr_i \Gamma$ we have $T = S \cap \pr_i^{-1}(\gamma)$.
\end{definition}

Now we can state the uniform Mordell--Lang theorem, which can be deduced from Theorem \ref{mordell-lang-combined} by automatic uniformity.

\begin{theorem}\label{uniform-Mordell-Lang}
Let $\mathfrak{S}$ be a semiabelian variety or $Y(1)^n$, and let $\Gamma \subseteq \mathfrak{S}$ be a structure of finite rank. Given a parametric family $(V_q)_{ q\in Q}$ of algebraic subvarieties of $\s$, there are a finite collection $\Sigma$ of special subvarieties of $\mathfrak{S}$ and an integer $m$, such that for every $q\in Q$ the variety $V_q$ contains at most $m$ maximal $\Gamma$-special subvarieties, each of which is a $\Gamma$-translate of a variety from $\Sigma$.
\end{theorem}
\begin{proof}
For the semiabelian case see \cite[Corollary 3.5.9]{hrushovski-manin-mumford} and \cite[Theorem 4.7]{scanlon-auto-uni}. For $Y(1)^n$ the theorem follows from Theorem \ref{mordell-lang-combined} and \cite[Theorem 2.4]{scanlon-auto-uni}, since $\Gamma$-special points are Zariski dense in $\Gamma$-special subvarieties.
\end{proof}

The following is a uniform version of our main theorems.

\begin{theorem}\label{gamma-cit-uniform}
Let $\mathfrak{S}$ be a semiabelian variety or $Y(1)^n$, let $\Gamma\seq \s$ be a structure of finite rank, and let $S \subseteq \mathfrak{S}$ be a $\Gamma$-special subvariety. Given a parametric family  $(V_q)_{ q\in Q}$ of algebraic subvarieties of $S$, there are a finite collection $\Sigma$ of special subvarieties of $S$ and an integer $m$, such that for any $q\in Q$ there is a finite subset $\Delta(q) \subseteq \Gamma$, with $|\Delta(q)| \leq m$, such that any $\Gamma$-atypical subvariety of $V_q$ is contained in a $\Delta(q)$-translate $T$ of some special variety from $\Sigma$ and $T\subsetneq S$.
\end{theorem}

\begin{proof}
We assume $\s$ is a semiabelian variety. The case of $Y(1)^n$ is completely analogous.

The proofs of Theorems \ref{gamma-cit} and \ref{gamma-semiab-special} can be generalised to work in this setting. In particular, we may assume $S = \mathfrak{S}$. Note that for a parametric family $(V_q)_{q\in Q}$, there is a parametric family consisting of all irreducible components of $V_q$ for all $q$,\footnote{By the results of \cite{vandenDries-bounds} the number of irreducible components of varieties in a parametric family is bounded, hence they form a parametric family as well.} hence we may assume each $V_q$ is irreducible. Further, let $T$ be one of the finitely many semiabelian subvarieties of $\mathfrak{S}$ given by \cite[Theorem 4.6]{Kirby-semiab}. Then the varieties $C_{T,V_q}$ defined as in Lemma \ref{cit-generic-typical} form a parametric family and we apply Theorem \ref{uniform-Mordell-Lang} to that family and proceed as in the proof of Theorem \ref{gamma-cit}. 
\end{proof}

\section{Optimal varieties}\label{section-optimal}

The Zilber--Pink conjecture is often formulated in terms of optimal subvarieties. Let $\mathfrak{S}$ be a semiabelian variety or $Y(1)^n$.

\begin{definition}[\cite{Pink,Habegger-Pila-o-min-certain}]
\begin{itemize}[leftmargin=0.5cm]
    \item[]
    \item For a subvariety $X\subseteq \mathfrak{S}$ the \emph{defect} of $X$ is the number $\delta (X) := \dim \langle X \rangle - \dim X.$
    \item Let $V$ be a subvariety of $\mathfrak{S}$. A subvariety $X \subseteq V$ is \emph{optimal} (in $V$) if for every subvariety $Y\subseteq V$ with $X \subsetneq Y$ we have $\delta(Y) > \delta(X)$.
\end{itemize}
\end{definition}

Observe that maximal atypical subvarieties are optimal, and optimal subvarieties are atypical but not necessarily maximal atypical. 

\begin{conjecture}[\cite{Habegger-Pila-o-min-certain}]\label{ZP-optimal}
Let $V$ be a subvariety of $\mathfrak{S}$. Then $V$ contains only finitely many optimal subvarieties.
\end{conjecture}

By \cite[Lemma 2.7]{Habegger-Pila-o-min-certain} this is equivalent to the Zilber--Pink conjecture. By analogy with optimal varieties, we want to define $\Gamma$-optimal varieties for a structure $\Gamma\seq \s$ of finite rank. For simplicity we focus on $Y(1)^n$.

\begin{definition}
Let $X$ be a subvariety of $Y(1)^n$.
\begin{itemize}[leftmargin=0.5cm]
    \item The $\Gamma$-\emph{special closure} of $X$, denoted $\langle X \rangle_{\Gamma}$, is the smallest $\Gamma$-special subvariety of $Y(1)^n$ containing $X$.
    \item The $\Gamma$-\emph{defect} of $X$ is the number $\delta_{\Gamma} (X) := \dim \langle X \rangle_{\Gamma} - \dim X.$
\end{itemize}

\begin{remark}
It is easy to verify that irreducible components of a non-empty intersection of $\Gamma$-special varieties are $\Gamma$-special, hence the $\Gamma$-special closure is well defined.
\end{remark}

\begin{definition}
Let $V$ be a subvariety of $Y(1)^n$ and $X$ be a subvariety of $V$. Then $X$ is called $\Gamma$-\emph{optimal} (in $V$) if whenever $X \subsetneq Y \subseteq V$, we have $\delta_{\Gamma} (X) < \delta_{\Gamma} (Y)$.
\end{definition}

\end{definition}

\begin{theorem}\label{thm-gamma-optimal}
Let $V \subseteq Y(1)^n$ be a subvariety. Then $V$ contains only finitely many $\Gamma$-optimal subvarieties whose weakly special closure is $\Gamma$-special.
\end{theorem}
\begin{proof}
The obvious adaptation of the proof of \cite[Lemma 2.7]{Habegger-Pila-o-min-certain} works in this setting.
\end{proof}

In the case of semiabelian varieties, and even algebraic tori, the irreducible components of an intersection of $\Gamma$-special subvarieties may not be $\Gamma$-special, hence we cannot define a $\Gamma$-special closure as above. Indeed, consider the two dimensional torus $\gm^2(\mathbb{C}) = (\mathbb{C}^{\times})^2$. Let $\Gamma_1$ be the torsion subgroup of $\gm(\mathbb{C})$, and let $\Gamma_2$ be the division closure of a cyclic subgroup of $\gm(\mathbb{C})$ generated by a transcendental element $\gamma \in \mathbb{C}$. Let also $\Gamma := \Gamma_1 \times \Gamma_2 \subseteq \gm^2(\mathbb{C})$. Consider two $\Gamma$-special subvarieties $ S: y_1y_2 = \gamma \mbox{ and } T: y_1^2y_2 = \gamma^2.$ Then $S \cap T = \{ (\gamma, 1) \}$ which does not contain a point of $\Gamma$, for $\gamma$ is not a torsion point.

However, in some cases the $\Gamma$-special closure is well-defined, and then the analogue of Theorem \ref{thm-gamma-optimal} clearly holds. For instance, when $\Gamma \subseteq \mathfrak{S}$ is the torsion subgroup of a semiabelian variety $\mathfrak{S}$, then $\Gamma$-special varieties coincide with special varieties and the $\Gamma$-special closure of an irreducible variety is equal to its special closure and is well-defined. In this case, the analogue of Theorem \ref{thm-gamma-optimal} states that for every variety $V\subseteq \mathfrak{S}$ there are only finitely many optimal subvarieties of $V$ whose weakly special closures are special (and one can use the Manin--Mumford conjecture instead of the Mordell--Lang conjecture to prove this). In the case of abelian varieties this is Corollary 9.11 of \cite{Habegger-Pila-o-min-certain}.

Let us give one more example when the $\Gamma$-special closure is well-defined. If $\mathfrak{S} = \gm^n$ is an $n$-dimensional torus, and $\Gamma = \Gamma_0^n$ where $\Gamma_0 \subseteq \gm$ is the division closure of a finitely generated subgroup (this is a direct analogue of a structure of finite rank in $Y(1)^n$), then it is easy to verify that $\Gamma$-special varieties are closed under taking irreducible components of intersections. Hence, the analogue of Theorem \ref{thm-gamma-optimal} holds in this case too.

As mentioned in the introduction, our methods are quite general and we expect them to extend to the setting of (pure) Shimura varieties, and the analogue of Theorem \ref{thm-gamma-optimal} should follow from an appropriate Ax--Schanuel statement (which was proven for pure Shimura varieties in \cite{Mok-Pila-Tsim}) and a Mordell--Lang conjecture (see, for example, \cite{Daw-Ren,Habegger-Pila-o-min-certain} for a discussion of the Zilber--Pink conjecture for Shimura varieties and the appropriate definitions in that setting). Further,  in \cite{Daw-Ren} Daw and Ren proved that the Zilber--Pink conjecture for Shimura varieties can be reduced to a point counting conjecture stating that every variety contains only finitely many optimal points. It seems their methods can be applied to prove an analogue of Theorem \ref{thm-gamma-optimal} for Shimura varieties.



I discussed these ideas with Christopher Daw, and he showed in particular that the argument of \cite[Theorem 8.3]{Daw-Ren} can be adapted to prove that if every variety contains only finitely many points which are special and optimal, then every variety contains only finitely many optimal subvarieties whose weakly special closures are special. On the other hand, finiteness of special optimal points follows from the Andr\'e--Oort conjecture for such points are maximal special. Thus, the Andr\'e--Oort conjecture for Shimura varieties implies that a subvariety of a Shimura variety contains only finitely many optimal subvarieties the weakly special closures of which are special. Since the Andr\'e--Oort conjecture is proven for $\mathcal{A}_g$ (see \cite{Tsimerman-AO-A_g}), this gives an unconditional result in that case. This method should probably extend to $\Gamma$-special and $\Gamma$-optimal varieties which will then give a new proof for Theorem \ref{thm-gamma-optimal}, and hence for Theorem \ref{gamma-cit-intro} too. Nevertheless, we do not consider these questions in this paper.

\subsection*{Acknowledgements} I would like to thank Christopher Daw and Sebastian Eterovi\'c for useful discussions and comments. I am also grateful to the referees for valuable comments that helped me improve the presentation of the paper.



\addcontentsline {toc} {section} {Bibliography}
\bibliographystyle {alpha}
\bibliography {ref}

\begin{thebibliography}{vdDS84}

\bibitem[Asl18]{Aslanyan-weakMZPD}
Vahagn Aslanyan.
\newblock Weak {M}odular {Z}ilber-{P}ink with {D}erivatives.
\newblock Preprint, arXiv:1803.05895, 2018.

\bibitem[BMZ07]{Bom-Mas-Zan}
Enrico Bombieri, David Masser, and Umberto Zannier.
\newblock Anomalous subvarieties - structure theorems and applications.
\newblock {\em IMRN}, 19, 2007.

\bibitem[DR18]{Daw-Ren}
Christopher Daw and Jinbo Ren.
\newblock Applications of the hyperbolic {A}x-{S}chanuel conjecture.
\newblock {\em Compositio Mathematica}, 154(9):1843--1888, 2018.

\bibitem[HP12]{Habegger-Pila-beyond}
Philipp Habegger and Jonathan Pila.
\newblock Some unlikely intersections beyond {A}ndr\'{e}-{O}ort.
\newblock {\em Compos. Math.}, 148(1):1--27, 2012.

\bibitem[HP16]{Habegger-Pila-o-min-certain}
Philipp Habegger and Jonathan Pila.
\newblock O-minimality and certain atypical intersections.
\newblock {\em Ann. Sci. \'{E}c. Norm. Sup\'{e}r.}, 49(4):813--858, 2016.

\bibitem[Hru01]{hrushovski-manin-mumford}
Ehud Hrushovski.
\newblock The {M}anin--{M}umford conjecture and the model theory of difference
  fields.
\newblock {\em Ann. Pure Appl. Logic}, 112(1):43--115, 2001.

\bibitem[Kir09]{Kirby-semiab}
Jonathan Kirby.
\newblock The theory of the exponential differential equations of semiabelian
  verieties.
\newblock {\em Selecta Mathematica}, 15(3):445--486, 2009.

\bibitem[Lan73]{Lang-elliptic}
Serge Lang.
\newblock {\em Elliptic functions}.
\newblock Addison-Wesley, 1973.

\bibitem[Mar02]{Mar}
David Marker.
\newblock {\em Model Theory: An Introduction}.
\newblock Springer, 2002.

\bibitem[McQ95]{mcquillan}
Michael McQuillan.
\newblock Division points on semi-abelian varieties.
\newblock {\em Invent. Math.}, 120(1):143--159, 1995.

\bibitem[MPT19]{Mok-Pila-Tsim}
Ngaiming Mok, Jonathan Pila, and Jacob Tsimerman.
\newblock Ax-{S}chanuel for {S}himura varieties.
\newblock {\em Annals of Mathematics}, 189(3):945--978, 2019.

\bibitem[Pil14]{pila-ellipt-mod-surf}
Jonathan Pila.
\newblock Special point problems with elliptic modular surfaces.
\newblock {\em Mathematika}, 60(1):1--31, 2014.

\bibitem[Pin05a]{Pink-2}
Richard Pink.
\newblock A combination of the conjectures of {M}ordell-{L}ang and
  {A}ndr\'e-{O}ort.
\newblock In F.~Bogomolov and Y.~Tschinkel, editors, {\em Geometric methods in
  algebra and number theory}, volume 235, pages 251--282. Progress in
  Mathematics, Birkh\"auser Boston, 2005.

\bibitem[Pin05b]{Pink}
Richard Pink.
\newblock A common generalization of the conjectures of {A}ndr\'e-{O}ort,
  {M}anin-{M}umford and {M}ordell-{L}ang.
\newblock Available at \url{https://people.math.ethz.ch/~pink/ftp/AOMMML.pdf},
  2005.

\bibitem[PT16]{Pila-Tsim-Ax-j}
Jonathan Pila and Jacob Tsimerman.
\newblock Ax-{S}chanuel for the $j$-function.
\newblock {\em Duke Math. J.}, 165(13):2587--2605, 2016.

\bibitem[Sca04]{scanlon-auto-uni}
Thomas Scanlon.
\newblock Automatic uniformity.
\newblock {\em IMRN}, (62):3317--3326, 2004.

\bibitem[Sca18]{Pila-Scanlon-dif-ZP}
Thomas Scanlon.
\newblock Differential algebraic {Z}ilber-{P}ink theorems.
\newblock Available at
  \url{math.berkeley.edu/~scanlon/papers/ZP-Oxford-July-2018.pdf}, 2018.

\bibitem[Sha13]{Shafarevich}
Igor Shafarevich.
\newblock {\em Basic Algebraic Geometry}, volume~1.
\newblock Springer, third edition, 2013.

\bibitem[Tsi18]{Tsimerman-AO-A_g}
Jacob Tsimerman.
\newblock The {A}ndr\'{e}-{O}ort conjecture for {$\mathcal{A}_g$}.
\newblock {\em Ann. of Math. (2)}, 187(2):379--390, 2018.

\bibitem[vdDS84]{vandenDries-bounds}
Lou van~den Dries and Karsten Schmidt.
\newblock Bounds in the theory of polynomial rings over fields. {A} nonstandard
  approach.
\newblock {\em Invent. Math.}, 76(1):77--91, 1984.

\bibitem[Zan12]{Zannier-book-unlikely}
Umberto Zannier.
\newblock {\em Some problems of unlikely intersections in arithmetic and
  geometry}, volume 181 of {\em Annals of Mathematics Studies}.
\newblock Princeton University Press, Princeton, NJ, 2012.
\newblock With appendixes by David Masser.

\bibitem[Zil02]{Zilb-exp-sum-published}
Boris Zilber.
\newblock Exponential sums equations and the {S}chanuel conjecture.
\newblock {\em J.L.M.S.}, 65(2):27--44, 2002.

\end{thebibliography}

\end{document}